\def\Bbb#1{{\bf #1}}

\def\fnote#1{\footnote}
\def\blacksquare{\hbox{\vrule width 4pt height 4pt depth 0pt}}

\def\cwleftpar#1#2{\leftskip #1 \rightskip #2 plus 1fill}
\def\cwrightpar#1#2{\leftskip #1 plus 1fill \rightskip #2}
\def\cwcenterpar#1#2{\leftskip #1 plus 1fill \rightskip #2 plus 1fill}
\def\cwfullpar#1#2{\leftskip#1\rightskip#2}

\def\cwoutdent#1#2{\llap{\hbox to #1{#2 \hss}}\ignorespaces}
\def\cwparbegin#1#2#3#4#5{
	\ifcase #1 \cwleftpar{#2}{#3}
	\or \cwrightpar{#2}{#3}
	\or \cwcenterpar{#2}{#3}
	\else \cwfullpar{#2}{#3}\fi
	\ifcase #4 \baselineskip = 1.5\baselineskip
	\or \baselineskip = 2\baselineskip
	\or \baselineskip = 3\baselineskip
	\else \baselineskip = 1\baselineskip\fi
	\ifdim #5 > 0in \else \noindent \fi
	\noindent\ignorespaces}
\documentclass{article}
\begin{document}
\advance \vsize by -1\baselineskip
\def\makefootline{
{\vskip \baselineskip \noindent \folio                                  \par
}}

\vspace*{2ex}
\noindent {\Huge Transports along Paths in\\[0.4ex] Fibre Bundles}\\[1.3ex]
\noindent {\Large I. General Theory}

\vspace*{2ex}

\noindent Bozhidar Zakhariev Iliev
\fnote{0}{\noindent $^{\hbox{}}$Permanent address:
Laboratory of Mathematical Modeling in Physics,
Institute for Nuclear Research and \mbox{Nuclear} Energy,
Bulgarian Academy of Sciences,
Boul.\ Tzarigradsko chauss\'ee~72, 1784 Sofia, Bulgaria\\
\indent E-mail address: bozho@inrne.bas.bg\\
\indent URL: http://theo.inrne.bas.bg/$\sim$bozho/}

\vspace*{2ex}

{\bf \noindent Published: Communication JINR, E5-93-299, Dubna, 1993}\\[1ex]
\hphantom{\bf Published: }
http://www.arXiv.org e-Print archive No.~math.DG/0503005\\[2ex]

2003 PACS numbers: 02.40.Ma, 04.90.+e\\[2ex]

\noindent
{\small
The \LaTeXe\ source file of this paper was produced by converting a
ChiWriter 3.16 source file into
ChiWriter 4.0 file and then converting the latter file into a
\LaTeX\ 2.09 source file, which was manually edited for correcting numerous
errors and for improving the appearance of the text.  As a result of this
procedure, some errors in the text may exist.
}\\[2ex]

	\begin{abstract}
Transports along path in fibre bundles are axiomatically introduced. Their
general functional form and some their simple properties are investigated.
The relationships of the transports along paths and lifting of paths are
studied.
	\end{abstract}\vspace{3ex}

{\bf 1. INTRODUCTION}
\nopagebreak

\medskip
The parallel transport (translation) is a well known concept in differential
geometry and fibre bundle theory, usually tied up with the connection theory
$[1-4]$. The purpose of this paper is to propose an independent definition
and the corresponding investigation of (parallel) transports (along paths) in
arbitrary fibre bundles. It is written in a manner analogous to the one of
[5,6], from where some ideas, results and proofs are transferred mutatis
mutandis. The detailed comparison of the material presented here with the one
in the literature on the same subject will be done in the next part of our
series.

The transports along paths in fibre bundles are defined in Sect. 2, where also two main groups of possible restrictions on them are considered. Some simple properties and the general functional form of these transports are found in Sect. 3. Sect. 4 deals with different connections between transports along paths and liftings of paths from the base to the bundle space of the fibre bundle. Necessary and sufficient condition is obtained for a fibre bundle to admit transports along paths. In Sect. $5, a$ method for generating a connection by a transport along a path is considered.

\medskip
\medskip
 {\bf 2. TRANSPORTS ALONG PATHS  AND POSSIBLE}

 {\bf RESTRICTIONS ON THEM}

\medskip
Let $(E,\pi ,B)$ be a general topological fibre bundle with a base $B$, total bundle space $E$ and projection $\pi :E  \to B [7-9]$. The fibre
bundle $(E,\pi ,B)$ is, generally, not supposed to be locally trivial. The fibres $\pi ^{-1}(x), x\in B$ are supposed to be homeomorphic with each other. The set of all sections of $(E,\pi ,B)$ is denoted by Sec$(E,\pi ,B)$, i.e., $\sigma \in $Sec$(E,\pi ,B)$ means $\sigma :B  \to E$ and $\pi \circ \sigma =${\it id}$_{B}[7,9]$, where {\it id}$_{X}$is the identity map of the set X.

By $J$ and $\gamma :J  \to B$ are denoted arbitrary, respectively, real interval and a path in B.

The transports along paths in $(E,\pi ,B)$ are defined in Subsect. 2.1 where also their basic properties, describing to a certain extend the dependence on their parameters, are found. On the transports along paths one can be imposed different restrictions. In the present work two main groups of them are considered. Firstly (Subsect. 2.2), the one describing more or less the transport's dependence on the path of the transport. The second one (Subsect. 2.3) concerns the ties between some primary given structures on the fibre bundle, if any, and transports along paths in it.

\medskip
\medskip
 {\bf 2.1. DEFINITION OF TRANSPORTS ALONG PATHS}

 {\bf IN FIBRE BUNDLES}

\medskip
The analysis of definition 2.1 of Ref. [6] shows that it may be generalized in such a way as the defined in it linear transport along paths in vector bundles to be a special case of more general "transports along paths" described by

{\bf Definition} ${\bf 2}{\bf .}{\bf 1}{\bf .} A$ transport along paths in
the fibre bundle $(E,\pi ,B)$ is a map I which to any path $\gamma :J  \to B$
puts into correspondence a map $I^{\gamma }$, transport along $\gamma $ such
that $I^{\gamma }:(s,t)\to I^{\gamma }_{s  \to t}$, where for every $s,t\in
J$ the map
\[
 I^{\gamma }_{s  \to t}:\pi ^{-1}(\gamma (s))  \to \pi ^{-1}(\gamma
(t)),\qquad (2.1)
\]
 a transport along $\gamma $ from $s$ to $t$, has the
following two properties:
\[
 I^{\gamma }_{t  \to r}\circ I^{\gamma }_{s  \to t}=I^{\gamma }_{s  \to r},
\quad r,s,t\in J,\qquad (2.2)
\]
\[
I^{\gamma }_{s  \to s}={\it id}_{\pi ^{-1_{(\gamma (s))}}},\quad s\in J.
 \qquad (2.3)
\]
 The path $\gamma $ and the numbers $s,t\in J$ in the map (2.1)
will be called, respectively, path, initial parameter and final parameter of
the transport.

 The property (2.2), which may be called a group property of transports along
paths, is an exact expression of the representation that the "composition of
two transports along one and the same path" must be a "transport along the
same path". The property (2.3) fixes a 0-ary operation in the set of
"transports along paths" defining in it the "unit transport" and, besides, it
is an exact expression of the naive understanding that if we "stand" at one
point of a path without "moving" along it, then "nothing must happen" with a
fibre over it.

The problem of existence of transports along paths, i.e. when a fibre bundle admits transports along paths, will be considered in Sect. 4 and the problem of their uniqueness or a general form will be investigated in Sect. 3.

Closely connected with definition 2.1 is the following one generalizing the concept for sections linearly transported along paths $(cf. [10]$, Sect. 5).

{\bf Definition 2.2.} Let in $(E,\pi ,B)$ be given a transport along paths
I. The section $\sigma \in $Sec$(E,\pi ,B)$ undergoes a (an I) transport or
it is (I-)transported (resp. along $\gamma :J\to B)$, if the equality
\[
\sigma (\gamma (t))=I^{\gamma }_{s\to t}\sigma (\gamma (s)),\quad s,t\in
J\qquad (2.4)
\]
 holds for every (resp. the given) path $\gamma :J\to $B.

{\bf Proposition 2.1.} If (2.1) is fulfilled for a fixed value $s\in J$,
then this equality is valid for every $s\in $J.

 {\bf Proof.} This result is a trivial corollary of (2.2).\blacksquare

Proposition 2.1 shows that in definition 2.2 it is sufficient to want (2.4) to be valid for a fixed $s=s_{0}\in J$ and then (2.4) may be regarded as a necessary and sufficient condition for the section $\sigma $ to be I-transported along $\gamma $.

{\bf Proposition 2.2.} If the section $\sigma \in $Sec$(E,\pi ,B)$ is I-transported along $\gamma :J\to B$, then its values $\sigma (\gamma (s))$ for every $s\in J$ are uniquely defined if the value $\sigma (\gamma (s_{0}))$ is given for an arbitrary fixed $s_{0}\in $J.

 {\bf Proof.} This result follows from (2.4) for $s=s_{0}.\blacksquare $

Evident example of transports along paths are the linear transport along paths in vector bundles, a result following from the comparison of definitions 2.1 from this work and from [10].

Now we shall consider two examples for transports along parts which appear in fibre bundles with a certain structure.

{\bf Example 2.1.} Let the fibre bundle $(E,\pi ,B)$ have a structure of
foliation [11], i.e. on the total bundle space $E$, which now is supposed to
be a manifold, to be given a foliation $\{K_{\alpha }: K_{\alpha }\subset E,
\alpha \in A\} [11]$, which, in particular, means that
$K_{\alpha }\cap K_{\beta }=\emptyset$,
$\alpha ,\beta \in A$, $\alpha \neq \beta $
 and $\cap_{\alpha \in A} K_{\alpha }=E$. Let the foliation $\{K_{\alpha
}\}$ be such that $\pi (K_{\alpha })=B, \alpha \in $A.

Before going on we want to stress the fact that this construction is equivalent to the one when over $(E,\pi ,B)$ is defined a family of sections $\{\sigma _{\alpha }: \sigma _{\alpha }\in $Sec$(E,\pi ,B), \alpha \in A\}$ such that if $\sigma _{\alpha }(x)=\sigma _{\beta }(x)$ for some $x\in B$, then $\alpha =\beta $. Actually, if such a family $\{\sigma _{\alpha }\}$ is given, then it is sufficient to put $K_{\alpha }=\sigma _{\alpha }(B), \alpha \in A$ and on the opposite, if a foliation $\{K_{\alpha }\}$ is given, then the sections $\sigma _{\alpha }, \alpha \in A$ are defined by the equality $\sigma _{\alpha }(x):=\pi ^{-1}(x)$  $K_{\alpha }, x\in B, \alpha \in $A.

In fibre bundles with such a structure there appears a natural concept for a lifting $\bar{\gamma }_{u}:J\to E$ of any path $\gamma :J\to B$ in the base $B$ through every point $u$ over the set $\gamma (J)$, i.e. $u\in \bar{\gamma }_{u}(J)$ and $\pi \circ \bar{\gamma }_{u}=\gamma
[12]$. In fact, if we define $\alpha (u)\in A$ as the subscript of the
unique fibre $K_{\alpha (u)}\in \{K_{\alpha }\}$ to which belongs the point
$u$, i.e. $K_{\alpha (u)}\ni u$, then the lifting $\bar{\gamma }_{u}:J\to E$
of $\gamma $ through $u$ is given by
\[
\bar{\gamma }_{u}(s)
:= \pi ^{-1}(\gamma (s)) \bigcap K_{\alpha (u)},\quad s\in J.
\]

Evidently, $\{\bar{\gamma }_{u}(J): u\in \pi ^{-1}(\gamma (s_{0}))\}$ for a
fixed $s_{0}\in J$ is a one dimensional foliation of the
$(\dim(B)+1)$-dimensional manifold $\pi ^{-1}(\gamma (J))$. Besides, it is
clear that $\bar{\gamma }_{u}$is the unique lifting of $\gamma $ in $E$ lying
as a whole in some of the fibres of the foliation $\{K_{\alpha }\}$.

The so-defined lifting $\gamma \to \bar{\gamma }_{u}$generates a transport
$K$ along the parts in $B$, defined by
\[
K^{\gamma }_{s\to t}(u):=\bar{\gamma }_{u}(t), u\in \pi ^{-1}(\gamma (s)),
\quad s\in J.
\]

The following equalities, which are true for every $r,s,t\in J$ and $u\in
\pi ^{-1}(\gamma (s))$, show that $K^{\gamma }_{s\to t}$is really a transport
along $\gamma $ from $s$ to $t ($see definition 2.1):
\[
K^{\gamma }_{s\to s}(u)=\bar{\gamma }_{u}(s)
=\pi ^{-1}(\gamma (s))
\bigcap K_{\alpha (u)}
=u,\quad  u\in \pi ^{-1}(\gamma (s)),\ s\in J,
\]
\[
 ={\bigl(}K^{\gamma }_{t\to r}\circ K^{\gamma }_{s\to t}{\bigr)}(u)
=K^{\gamma }_{t\to r}{\bigl(}\bar{\gamma }_{u}(t){\bigr)}
=\bar{\gamma }_{\bar{\gamma}_u(t)}(r)
=\pi ^{-1}(\gamma (r)) \bigcap K_{\alpha (\bar{\gamma}_u(t))}
\]
\[
=\pi ^{-1}(\gamma (r)) \bigcap K_{\alpha (u)}
=\bar{\gamma }_{u}(r)
=K^{\gamma }_{s\to r}(u),
\quad  u\in \pi ^{-1}(\gamma (s)),\ r,s,t\in J.
\]

 Here, we have used that $\alpha (v)=\alpha (u)$ for every $v\in K_{\alpha
(u)}$which is a consequence of the uniqueness of the fibre  of the foliation
passing through an arbitrary point of the foliation.

{\bf Example 2.2.} In [1], vol I, p. 174 the  definitions are given for
parallelism and a parallel vector field according to which  "if $\xi $ is a
vector bundle over $B\times B$, whose bundle over $(x,y)\in B\times B$ is
$\{f:  f:T_{x}(B)\to T_{y}(B), f$ - linear$\}$, then the parallelism on $B$
is  such a section $P\in $Sec$\xi $ that \[ P(z,y)\circ P(x,z)=P(x,y),\quad
P(x,x)={\it id}_{T_{x}},\ x,y,z\in B.  \] The vector field $X$ is parallel
(with respect to $P)$ if \[ P(x,y)X_{x}=X_{y}".  \]

If $B$ is parallelizable, i.e. in it exists a parallelization $P$, then it uniquely defines axiomatically defined (global) parallel transport in $B [2]$. From the cited definitions it is clear that the parallelization $P$ on $B$ defines a transport (along paths) in the tangent to $B$ fibre bundle $(T(B),\pi ,B)$. Besides, this transport is global in a sense that it depends only on the initial and final points of a transport and does not depend on paths connecting them. As the vector fields over $B$ are in fact sections of $(T(B),\pi ,B) [1,2]$, it is evident that the above definition for parallelism of $X$ is a special case of the definition 2.2.

\medskip
\medskip
 {\bf 2.2. ADDITIONAL RESTRICTIONS}
\nopagebreak

\medskip
The considered below restrictions, which can be imposed on the transports along paths, are under the title "additional" as, nevertheless they are important from some view-points, in our opinion, they have more special character and define one or another special property of the transports along paths. One should have in mind that they describe the functional dependence of the transports along paths on the path or curve of transport and they can be imposed separately, on groups, as well as all of them together. Besides, they have a number of equivalent formulations, but we below consider only one of them.

 {\bf (1)} {\sl Condition for locality}

 If $J^\prime $ is a subinterval of $J$ and $\gamma :J\to B$, then
 \[
 I^{\gamma \mid J^\prime }_{s\to t}=I^{\gamma }_{s\to t}. s,t\in J^\prime
\qquad (2.5)
\]
 where $\gamma \mid J^\prime $ is the restriction of $\gamma $
on $J^\prime \subset $J.

 {\bf (2)} {\sl Condition for invariance under the parameter changes}

 If $\tau :J^{\prime\prime}\to J$ is one-to-one map and $\gamma :J\to B$,
then
\[
I^{\gamma o\tau }_{s\to t}=I^{\gamma }_{\tau (s)\to \tau (t)}, s,t\in
J^{\prime\prime}.\qquad (2.6)
\]

 {\bf (}{\bf 3}{\bf )} {\sl Conditions for smoothness}

We will not formulate these conditions, but we shall only mention that there are three types of them, describing the "smoothness" respectively of:

 {\bf (3a)} the dependence of a transport on a path of the transport;

{\bf (3b)} the dependence of a transport on its initial and final parameters;

{\bf (3c)} a transport as a map between the fibres over the curve
defined by the path of transport.

In connection with the locality condition the following proposition is very important.

{\bf Proposition 2.3.} If $I^{\gamma }$is a transport along $\gamma $, then
(2.5) is equivalent to
\[
 I^{\gamma }_{s\to t}=I^{\gamma \mid [\min(s,t),\max(s,t)]}_{s\to t}, \quad
s,t\in J.  \qquad (2.7)
\]

 {\bf Proof.} (2.7) follows from (2.5) for
$J^\prime =[\min(s,t),\max(s,t)]\subset $J. On the opposite, if (2.7) holds,
then
\[
 I^{\gamma \mid J^\prime }_{s\to t}=I^{(\gamma \mid J^\prime )\mid
[\min(s,t),\max(s,t)]}_{s\to t}=I^{\gamma \mid [\min(s,t),\max(s,t)]}_{s\to
t}=I^{\gamma }_{s\to t},
\]
 where $s,t\in J^\prime \subset J$ and we have used the evident fact that
if  $J^{\prime\prime}\subset J^\prime \subset J$, then $(\gamma \mid J^\prime
)\mid J^{\prime\prime}=\gamma \mid J^{\prime\prime}.\blacksquare $

par
Proposition 2.3 shows that if (2.5) holds, then the transport along $\gamma $ from $s$ to $t$ does not depend globally on the "whole" map $\gamma $, but only on its restriction on the interval defined from $s$ and t.

The transports along paths satisfying the condition (2.6), in fact, depend not on the path of transport $\gamma :J\to B$, bun on the curve
of transport, i.e. only on the values $\gamma (s), s\in J$, or other wise stated, (2.6) leads to the dependence of a transport along $\gamma $ only on the nonparametrized curve $\gamma (J)$. Namely, the set of these transports has invariant and clear geometrical sense.

It is important to note that frequently the condition (2.6) is given in a weaker form in which one wants $\tau $ to be orientation preserving homeo- or diffeomorphism.

Elsewhere it will be proved that "a transport along paths is a parallel transport iff it satisfies simultaneously (2.5) and (2.6)".

The conditions of the type (3a) in their essence coincide with the conditions for smoothness in the axiomatic approach to the parallel transport, and those of the type (3b) in the case of considered in [6] linear transports in vector bundles are reduced to the standard condition for smoothness of linear operators (between vector spaces) on scalar parameter.

\medskip
\medskip
 {\bf 2.3. CONDITIONS FOR CONSISTENCY}

\medskip
The conditions for consistency describe connections between transports along paths in fibre bundles and, if there exist, specific for the fibre bundles structures. In the most general case of arbitrary fibre bundles there are no such structures due to which there are no conditions for consistency. They appear in more special cases, two of which are considered below. The general case of this problem will be a subject of another work.

Let $I^{\gamma }_{s\to t}$be a transport along $\gamma :J\to B$ from $s$ to $t, s,t\in J$ in the real (resp. complex) vector bundle $(E,\pi ,B)$. Then, on the transport I there can be imposed the following restriction.

 {\bf (4)} {\sl Condition for consistency with a vector structure}

 If $\lambda ,\mu \in {\Bbb R} ($resp. $\lambda ,\mu \in {\Bbb C})$ and
$u,v\in \pi ^{-1}(\gamma (s)), s\in J$, then
\[
 I^{\gamma }_{s\to t}(\lambda u+\mu v)=\lambda I^{\gamma }_{s\to t}u+\mu
I^{\gamma }_{s\to t}v,\quad s,t\in J.
 \qquad (2.8)
\]
 Using
definition 2.2 we can reformulate (2.8) in an equivalent way by saying that a
real (resp.  complex) linear combination of transported along $\gamma $
sections is also a transported along $\gamma $ section whose value at
$\gamma(t)$ is obtained by transporting along $\gamma $ from $s$ to $t$ the
value of this linear combination at a point $\gamma (s)$ for an arbitrary
fixed $s\in $J.

Evidently, in vector bundles the condition (2.8) reduces the transports along paths to the considered in [6] linear transports along paths.

Another important example of conditions for consistency appears in the case when in the fibre bundle $(E,\pi ,B)$ is given a bundle metric $g$, i.e. [3] when in it is fixed a family $\{g_{x}: g_{x}:\pi ^{-1}(x)  \pi ^{-1}(x)\to {\Bbb R}, x\in B\}$. In particular, if $(E,\pi ,B)$ is a vector fibre bundle and $g_{x}, x\in B$ are nondegenerate bilinear forms, then $g$ defines an usual bundle metric in the fibre bundle [2]. (To this class belong also all (pseudo-)Riemannian metrics: if $B$ is a (pseudo-)Riemannian manifold, then one puts $E=T(B), \pi ^{-1}(x)=T_{x}(B), x\in B$ and $g$ is identified with the (pseudo-)Riemannian metric of B.) In this case the class of transports along $\gamma :J\to B$ satisfying the following condition is naturally generated.

 {\bf (5)} {\sl Condition for consistency with a bundle metric}

 If $u,v\in \pi ^{-1}(\gamma (s)), s\in J$, then
\[
g_{\gamma (s)}(u,v)=g_{\gamma (t)}(I^{\gamma }_{s\to t}u,I^{\gamma }_{s\to
t}v),\quad s,t\in J.  \qquad (2.9)
\]
 In the (pseudo-)Riemannian case
(2.9) simply means a conservation of the scalar product of the vectors under
their (parallel) transport along an arbitrary path at any its point.

The transports along paths satisfying simultaneously (2.8) and (2.9) will be investigated elsewhere $(cf. [13])$.

As we shall see elsewhere the conditions for consistency of transports along paths with some structures are equivalent to the demand these structures to be "transported" with respect to properly chosen transports along paths acting in the corresponding fibre bundles. In particular, this is the case with the condition for consistency of a parallelism $P ($see example 2.2 at the end of Subsect. 2.1) and its torsion $S$, which is given in [1], vol. I, p.175 and p.235, where this condition is taken as a definition for parallelness ("parallel transport") of $S$; in fact, in this case, these two concepts are equivalent.

\medskip
\medskip
 {\bf 3. SOME PROPERTIES AND GENERAL FORM}

\medskip
 {\bf Proposition 3.1.} If (2.2) is valid, then (2.3) is equivalent to the
invertability of the map (2.1) when its initial and final parameters
coincide, i.e. to the existence of ${\bigl(}I^{\gamma }_{s\to s}
\bigr)^{-1}, s\in $J.

{\bf Proof.} If (2.3) is valid, then, evidently, no matter whether (2.2) is
true or not ${\bigl(}I^{\gamma }_{s\to s}{\bigr)}^{-1}$ there exists and
${\bigl(}I^{\gamma }_{s\to s}{\bigr)}^{-1}=${\it id}$_{\pi
^{-1_{(\gamma (s))}}}=I^{\gamma }_{s\to s}$. On the opposite, if, for
example, exists the left inverse map ${\bigl(}I^{\gamma }_{s\to s}
\bigr)^{-1}$, then putting $r=s=t$ in (2.2), we get $I^{\gamma }_{s\to
s}\circ I^{\gamma }_{s\to s}=I^{\gamma }_{s\to s}$ and multiplying from left
by ${\bigl(}I^{\gamma }_{s\to s}{\bigr)}^{-1}$, we convince ourselves in the
validity of (2.3).\blacksquare

 {\bf Proposition 3.2.} If (2.2) is fulfilled, then (2.3) is equivalent to
the existence of the inverse map of (2.1), besides
\[
 {\bigl(}I^{\gamma }_{s\to t}{\bigr)}^{-1}=I^{\gamma }_{t\to s},
\quad s,t\in J. \qquad (3.1)
\]

 {\bf Proof.} If (2.3)
holds, putting $r=t$ in (2.2), we get $I^{\gamma }_{s\to t}\circ \circ
I^{\gamma }_{t\to s}=${\it id}$_{\pi ^{-1_{(\gamma (t))}}}$, which results in
(3.1). On the opposite, if (3.1) is true, then putting $r=t$ in (2.2), we get
$I^{\gamma }_{t\to t}=I^{\gamma }_{s\to t}\circ  \circ I^{\gamma }_{t\to
s}=I^{\gamma }_{s\to t}\circ {\bigl(}I^{\gamma }_{s\to t} \bigr)^{-1}
=$id$_{\pi ^{-1_{(\gamma (t))}}}$, i.e. (2.3) is valid.\blacksquare

In particular, (3.1) shows that the transports along a path $\gamma :J\to B$ are 1:1 invertible maps between the fibres over $\gamma (J)$.

If $\tau _{-}:J^\prime \to J$ is an orientation changing 1:1 map between the real intervals $J^\prime $ and $J$ and $\gamma :J\to B$ is a  path, then the path $\gamma _{\tau _{-}}:=\gamma \circ \tau _{-}:J^\prime \to B$ will be called inverse to $\gamma $ with respect to $\tau _{-}$. In the canonical case, when paths of the type $\gamma :[0,1]  \to B$ are considered, by definition [8,12] the inverse  path to $\gamma :[0,1]:\to B$ is $\gamma _{-}:=\gamma \circ \tau ^{c}_{-}:[0,1]\to B$, where $\tau ^{c}_{-}:[0,1]\to [0,1]$ is defined by $\tau ^{c}_{-}:s\to 1-s, s\in [0,1]$.

{\bf Proposition 3.3.} If (2.6) holds, then the transports along a path
$\gamma :J\to B$ and its inverse path $\gamma _{\tau _{-}}:J^\prime \to B$
are connected by
\[
I^{\gamma _{\tau _{-}}}_{s\to t}=I^{\gamma }_{\tau _{-}},
\quad s,t\in J^\prime . \qquad (3.2)
\]

 {\bf Proof.} (3.2) follows from (2.6) for $\tau =\tau_{-}.\blacksquare $

We define the product of  paths in the following way which differs from the
standard one $(cf. [8,12]$; see bellow). Let us have an ordered pair $(\gamma
_{1},\gamma _{2})$ of two  paths $\gamma _{h}:[a_{h},b_{h}]\to B, h=1,2$ with
the end of $\gamma _{1}$coinciding with the beginning of $\gamma _{2}$, i.e.
$\gamma _{1}(b_{1})= =\gamma _{2}(a_{2})$. Let there be given numbers
$a_{o},b_{o},c_{o}\in {\Bbb R}$, such that $a_{o}\le c_{o}\le b_{o}$, and
one-to-one maps $\tau _{1}:[a_{o},c_{o}]\to [a_{1},b_{1}]$ and $\tau
_{2}:[c_{o},b_{o}]\to [a_{2},b_{2}]$ preserving the orientation, which, in
particular, means that $\tau _{1}(a_{o})=a_{1}, \tau _{1}(c_{o})=b_{1}, \tau
_{2}(c_{o})=a_{2}$and $\tau _{2}(b_{o})=b_{2}$. The product of the paths
$\gamma _{1}$and $\gamma _{2}$ is a  path
\[
 (\gamma _{1}\gamma _{2})_{\chi }:[a_{o},b_{o}]\to B,\qquad (3.3a)
\]
 which depends on the parameter $\chi :=(a_{o},b_{o},c_{o};\tau _{1},\tau
_{2})$, and it is defined by the equalities
\[
 (\gamma _{1}\gamma _{2})_{\chi } | [a_{o},c_{o}]:=\gamma _{1}\circ \tau
_{1},\qquad (3.3b)
\]
\[
 (\gamma _{1}\gamma _{2})_{\chi } | [c_{o},b_{o}]
:=\gamma _{2}\circ \tau _{2}.\qquad (3.3c)
\]
 In the canonical case (see [8,12]) the product of the
paths $\gamma _{h}:[0,1]\to B, h=1,2$ is defined by a parameter $\chi
^{c}=(0,1,1/2;\tau ^{c}_{1},\tau ^{c}_{2})$, where $\tau ^{c}_{1}:s\to 2s$
for $s\in [0,1/2]$ and $\tau ^{c}_{2}:s\to 2s-1$ for $s\in [1/2,1]$.

{\bf Proposition 3.4.} Let $(2.2), (2.6)$ and (2.7) be valid. Let $\gamma
_{h}:[a_{h},b_{h}]\to B, h=1,2$ and $(\gamma _{1}\gamma _{2})_{\chi }$be the
product of $\gamma _{1}$and $\gamma _{2}$defined by a parameter $\chi
=(a_{0},b_{0},c_{0};\tau _{1},\tau _{2})$. Then, for arbitrary $t_{1}\in
[a_{0},c_{0}]$ and $t_{2}\in [c_{0},b_{0}]$ the equality is valid
\[
I^{(\gamma _{1}}_{t_{1}\to t_2}
= I^{\gamma _{2}}_{a_{2}\to \tau_2(t_2)}\circ
  I^{\gamma _{1}}_{\tau_{1}(t_1)\to b_1},
\qquad  a_{0}\le t_{1}\le c_{0}\le t_{2}\le
b_{0}.  \qquad (3.4)
\]

 {\bf Proof.} Using consequently $(2.2), (2.7), (2.6)$ and the above
definition of $(\gamma _{1}\gamma _{2})_{\chi }$, we get:
\[
 I^{(\gamma _{1}\gamma_2)_\chi}_{t_{1}\to t_2}
= I^{(\gamma _{1}\gamma_2)_\chi}_{c_{0}\to t_2}\circ
  I^{(\gamma _{1}\gamma_2)_\chi}_{t_{1}\to c_0}
=I^{(\gamma _{1}\gamma_2)_\chi|[c_0,t_2]}_{c_{0}\to t_2}\circ
I^{(\gamma _{1}\gamma_2)_\chi|[t_1,c_0]}_{t_{1}\to c_0}
\]
\[
 = I^{\gamma _{2}\circ\tau_2}_{c_{0}\to t_2}\circ
I^{\gamma _{1}\circ\tau_1}_{t_{1}\to c_0}
= I^{\gamma_{2}}_{\tau_{2}(c_0)\to \tau_2(t_2)} \circ
  I^{\gamma_{1}}_{\tau_{1}(t_1)\to \tau_1(c_0)},
\]
 from where it follows (3.4), as we have $\tau _{1}(a_{0})=a_{1}, \tau
_{1}(c_{0})=b_{1}, \tau _{2}(c_{0})=a_{2}$and $\tau _{2}(b_{0})=b_{2}$because
$\tau _{1}$and $\tau _{2}$are orientation preserving .\blacksquare

{ \bf Remark.} If we had, for example, $t_{1},t_{2}\in [a_{0},c_{0}]$, then
instead of (3.4) the equality will be valid
\[
 I^{(\gamma _{1}\gamma_2)_\chi}_{t_{1}\to t_2}
= I^{\gamma _{1}}_{\tau _{1}\to\tau_1(t_2)},\quad  t_{1},t_{2}\in
[a_{0},c_{0}],\qquad (3.5)
\]
 in the derivation of which only (2.6) and (2.7) can be used.

The above proposition shows some important properties of transports along paths and connections between the basic properties describing them. Let us note that when deriving $(3.2), (3.4)$ and (3.5) the main role was played by the condition for invariance
(2.6).

Once again we want to stress the fact that the basic properties describe the dependence of a transport along paths on the initial and final value of its parameters, while the additional conditions define its dependence on the path of the transport.

The general functional form of the transports along paths is given by

{\bf Theorem 3.1.} Let in the base $B$ of the fibre bundle $(E,\pi ,B)$ be
given a path $\gamma :J\to B$ and for arbitrary $s,t\in J$ be given a map
(2.1). The maps $I^{\gamma }_{s\to t}, s,t\in J$ define a transport along
$\gamma $ from $s$ to $t$, i.e. (2.2) and (2.3) are satisfied iff there exist
a set $Q$ and a family of one-to-one maps $\{F^{\gamma }_{s}:\pi ^{-1}(\gamma
(s))\to Q, s\in J\}$ such that
\[
I^{\gamma }_{s\to t}={\bigl(}F^{\gamma }_{t}{\bigr)}^{-1}\circ
\bigl( F^{\gamma }_{s}{\bigr)},\quad s,t\in J. \qquad (3.6)
\]

{\bf Proof.} The theorem is a corollary of the following lemma in which one
has to put $N=J, Q_{s}=\pi ^{-1}(\gamma (s))$ and $R_{s  \to t}=I^{\gamma
}_{s\to t}, s,t\in $J.

{\bf Lemma 3.1.} Let there be given a set $N, N\neq\emptyset $ and families
of equipollent sets $\{Q_{s}: s\in N\}$ and of maps $\{R_{s\to t}: R_{s\to
t}:Q_{s}\to Q_{t}, s,t\in N\}$. Then, the maps of $\{R_{s\to t}\}$ satisfy
the equalities
\[
 R_{s\to t}\circ R_{r\to s}=R_{r\to t},\quad  r,s,t\in N,\qquad (3.7)
\]
\[
R_{s\to s}={\it id}_{Q_{s}},\quad  s\in N\qquad (3.8)
\]
 iff there exists an equipollent with $Q_{s}$for some $s\in N$ set $Q$ and a
family of one-to-one maps $\{F_{s}: F_{s}:Q_{s}\to Q, s\in N\}$, such that
\[
  R_{s\to t}=(F_{t})^{-1}\circ (F_{s}),\quad  s,t\in N.\blacksquare
 \qquad (3.9)
\]
 {\bf Proof.} The sufficiency is almost evident: the substitution of
(3.9) into (3.7) and (3.8) converts them into identities. On the opposite,
putting in $(3.7) r=t$ and using (3.8), we see that $R_{s\to t}$has an
inverse map and
\[
 (R_{s\to t})^{-1}=R_{t\to s},\quad  s,t\in N,\qquad (3.10)
\]
 due to which (see and
(3.7)) for any fixed $s_{0}\in N$, we have $R_{s\to t}= =R_{s_{0}}\circ
R_{s\to s_{0}}={\bigl(}R_{t\to s_{0}}{\bigr)}^{-1}\circ {\bigl(}R_{s\to
s_{0}}{\bigr)}$, i.e. (3.9) is fulfilled for $Q=Q_{s_{0}}$and
$F_{s}=R_{s\to s_{0}}.\blacksquare $

The arbitrariness in the choice of the set $Q$ and a family $\{F_{s}\}$ in the theorem 3.1 is described by

{\bf Proposition 3.5.} Let in the fibre bundle $(E,\pi ,B)$ be given a
transport along paths I with a representation (3.6) for some set $Q$ and a
family of 1:1 maps $\{F^{\gamma }_{s}:\pi ^{-1}(\gamma (s))\to Q, s\in J\}$.
Then, there exist a set $^{o}Q$ and a family of one-to-one maps
$\{^{o}F^{\gamma }_{s}:\pi ^{-1}(\gamma (s))\to ^{o}Q, s\in J\}$ such that
\[
 I^{\gamma }_{s\to t}={\bigl(}^{o}F^{\gamma }_{t}{\bigr)}^{-1}\circ
{\bigl(}^{o}F^{\gamma }_{s}{\bigr)}, s,t\in J,\qquad (3.6^\prime )
\]
 iff there exists an one-to-one map $D^{\gamma }:^{o}Q\to Q$ for which
\[
 F^{\gamma }_{s}=D^{\gamma }\circ {\bigl(}^{o}F^{\gamma }_{s}{\bigr)},
\quad s\in J. \qquad (3.11)
\]

 {\bf Proof.} The
proposition follows for $N=J, Q_{s}=\pi ^{-1}(\gamma (s))$ and $R_{s  \to
t}=I^{\gamma }_{s\to t}, s,t\in J$ from

{\bf Lemma 3.2.} Let there be given a set $N, N\neq   $, family of
equipollent sets $\{Q_{s}: s\in N\}$, an equipollent with $Q_{s}$for some
$s\in N$ set $^{o}Q$ and a family of maps $\{^{o}F_{s}:{ } ^{o}F_{s}:Q_{s}\to
^{o}Q, s\in N\}$. If
\[
  R_{s\to t}={\bigl(}^{o}F_{t}{\bigr)}^{-1}\circ
\bigl(^{o}F_{s}{\bigr)},\qquad (3.9^\prime )
\]
 then (3.9) is valid for some
family of maps $\{F_{s}: F_{s}:Q_{s}\to Q, s\in N\}, Q$ being some
equipollent with $^{o}Q$ set, if and only if there exists a one-to-one map
$D:^{o}Q\to Q$, such that
\[
F_{s}=D\circ (^{o}F_{s}).\blacksquare \qquad (3.12)
\]

 {\bf Proof.} The
sufficiency is almost evident: if (3.12) is true for some $Q$ and
$\{F_{s}\}$, then from it we find $^{o}F_{s}=D^{-1}\circ F_{s}, s\in N$ and
substituting this result into $(3.9^\prime )$, we get (3.9). On the opposite,
if (3.9) is true, then, due to $(3.9^\prime )$, from it results in
$\bigl( F_{t}{\bigr)}^{-1}\circ  \circ {\bigl(}F_{s}{\bigr)}=\bigl(^{o}
F_{t}{\bigr)}^{-1}\circ {\bigl(}^{o}F_{s}{\bigr)}, s,t\in N$,
hereof we see that $F_{s}\circ {\bigl(}^{o}F_{s}{\bigr)}^{-1}=
=F_{t}\circ {\bigl(}^{o}F_{t}{\bigr)}^{-1}$for any $s,t\in N$, but this
means that the left and right hand sides of the last equality do not depend
either on $s$ or on t. Hence, fixing arbitrarily some $s_{0}\in N$ and
putting $D=F_{s_{0}}\circ  \circ {\bigl(}^{o}F_{s_{0}}
\bigr)^{-1}:^{o}Q\to Q$, from the last equality for $t=s_{0}$, we get
(3.12).\blacksquare

\medskip
\medskip
 {\bf 4. TIES WITH THE LIFTING OF PATHS}

\medskip
The definition of a transport along paths I in a fibre bundle $(E,\pi ,B)$ leads to a natural lifting of every path $\gamma :J\to B [12]$.

{\bf Definition 4.1.} Let $\gamma :J\to B be a path, I^{\gamma }$ be a
transport along $\gamma$, $u\in\pi^{-1}(\gamma (J))$, $q_{\gamma }(u):=\{s:
s\in J$, $\gamma (s)=\pi (u)\}$ and $s_{0}\in q_{\gamma }(u)$. The path
$\bar{\gamma}_{u,s_{0}}:J\to E$ defined by the equality

\[
 \bar{\gamma}_{u,s_{0}}(s):=I^{\gamma }_{s_{0}\to s}(u),
\quad s\in J,\qquad (4.1)
\]
 is
the lifting of $\gamma $ through $u$ generated by $I^{\gamma }$ with a
parameter $s_{0}$; the map $\gamma \to \bar{\gamma}_{u,s_{0}}$ is the lift
through $u$ generated by $I^{\gamma }$ with a parameter $s_{0}$. If $q_{\gamma
}(u)$ consists of only one element, i.e., if $\pi (u)$ is not a
self-intersection point for $\gamma $, then   $\bar{\gamma}_{u,s_{0}}$ will
be denoted simply by   $\bar{\gamma}_{u}$ and we say that $\bar{\gamma}_{u}$
and $\gamma \to\bar{\gamma}_{u}$ are generated by $I^{\gamma }$.

Evidently, any transport along paths I generates through $(4.1) a$ lift of
the paths from the base B. The usage of the term "lifting" here and in
definition 4.1 is correct due to

 {\bf Proposition 4.1.} The path   $\bar{\gamma}_{u,s_{0}}:J\to E$, defined
by (4.1), is a lifting of the path $\gamma :J\to B$ (from $B$ into $E)$
through $u$, i.e.
\[ \pi \circ \bar{\gamma}_{u,s_{0}} = \bar{\gamma}_{u,s_{0}}(s_{0})=u  \qquad
(4.2)
\]

{\bf Proof.} (4.2) follows from (4.1) and definition 2.1 (see also (2.1) and
(2.3)).\blacksquare

From (4.1) we can make the following simple, but important, conclusion. In the general case, one transport along paths generates through the point $u\in \pi ^{-1}(\gamma (J))$ as many liftings of the path $\gamma :J\to B$ as is the number of elements (the power) of $q_{\gamma }(u)$, equal, evidently, to one plus the number of self-intersections of $\gamma $ at the point $\pi (u)$. In this connection, one naturally puts the question in what sense and when the lift of a given path through a lying above it point is unique.

{\bf Definition 4.2.} The lifting generated through (4.1) by a transport
along paths I is globally unique (or is unique in an absolute sense) if for
every path $\gamma :J\to B$ and every $u\in \pi ^{-1}(\gamma (J))$ we have
\[
\bar{\gamma}_{u,r}=\bar{\gamma}_{u,s},\quad r,s\in q_{\gamma }(u).\qquad
(4.3)
\]

{\bf Proposition 4.2.} The lifting generated through (4.1) by the transport
I is globally unique iff
\[
I^{\gamma }_{r\to s}={\it id}_{\pi ^{-1_{(\gamma (s))}}}
\quad for\ those\ r,s\in J
\ for\ which\ \gamma (r)=\gamma (s),  \qquad (4.4)
\]
 or, which is equivalent, iff
\[
v\in \bar{\gamma}_{u,r}(J)
 \Leftarrow \Leftrightarrow \bar{\gamma}_{v,s}= \bar{\gamma}_{u,r},
\quad s\in q_{\gamma }(v),\ r\in q_{\gamma }(u).\qquad (4.5)
\]

 {\bf Proof.} Let (4.3) be valid. Then, for any $t\in J$ from (4.1), we get
$I^{\gamma }_{r\to t}(u)=\bar{\gamma}_{u,r}(t)= \bar{\gamma}_{u,s}(t)
=I^{\gamma }_{s\to t}(u)$.  Therefore, due to the
arbitrariness of $u\in \pi ^{-1}(\gamma (J))$, we find $I^{\gamma }_{r\to
t}=I^{\gamma }_{s\to t}$, which as a consequence of $(2.1), (2.2), (3.1)$ and
the definition of $q_{\gamma }(u)$ (see definition 4.1), is equivalent to
(4.4). Further, if $v\in\bar{\gamma}_{u,r}(J)$, $r\in q_{\gamma }(u)$, then
by (4.1) there exists $s_{1}\in J$ such that
$v = \bar{\gamma}_{u,r}(s_{1})=I^{\gamma }_{r\to s_{1}}(u)$. So for every
$t\in J$ and $s\in q_{\gamma }(v)$, we have
\(
\bar{\gamma}_{v,s}(t)
=I^{\gamma }_{s\to t}(v)
=I^{\gamma }_{s\to t}\circ I^{\gamma }_{r\to s_{1}}(u)
=(I^{\gamma }_{s\to t}\circ
  I^{\gamma }_{s_{1}\to s}\circ I^{\gamma }_{r\to s_{1}})(u)
=I^{\gamma }_{r\to t}(u)=\bar{\gamma}_{u,r}(t).
\)
 Here we used the evident fact that $s_{1}\in q_{\gamma }(v)$ (see (4.2)),
the proved equality (4.4), and, besides, we have twice applied (2.2). Hence
$\bar{\gamma}_{v,s}=\bar{\gamma}_{u,r}$, i.e. the right implication in (4.5)
is true. The inverse implication in (4.5) is evident and it does not depend
on any additional facts: if $\bar{\gamma}_{v,s}=\bar{\gamma}_{u,r}$ for $s\in
q_{\gamma }(v)$ and $r\in q_{\gamma }(u)$, then
$v\in\bar{\gamma}_{v,s}(J)=\bar{\gamma}_{u,r}(J)$, i.e. $v\in $
$\bar{\gamma}_{u,r}(J)$.  So, we proved that (4.3) results in (4.4) leading
to (4.5) due to which the proof of the proposition ends with the fact that
from (4.5) there follows (4.3): as by definition $u\in \bar{\gamma}_{u,r}(J)$, then from (4.5) for $v=u$ there follows
(4.3).\blacksquare

In a local sense, the generated from a transport along paths lifting of paths may be unique under different criteria, modifying in an appropriate way definition 4.2 and to which there correspond modified versions of proposition 4.2. For instance, the local uniqueness may be defined in the following two ways:

(a) Equality (4.3) is valid only along a given path $\gamma  ($uniqueness along a given path). A necessary and sufficient condition for this is (4.4), or (4.5), to be fulfilled along the path $\gamma . ($The proof is a consequence form the one of proposition 4.2 for a fixed path $\gamma .)$

(b) Equality (4.3) to be valid only at a given point $u$ over a path $\gamma  ($uniqueness along a given path through a point above it). A necessary and sufficient condition for this is (4.4) to be fulfilled for every $r,s\in q_{\gamma }(u). ($The proof is a consequence from the one of proposition 4.2 for fixed $\gamma $ and u.) Due to this if the lifting is unique along $\gamma $ through $u$, then it is unique along $\gamma $ through every $v\in \pi ^{-1}(\pi (u))$, so $q_{\gamma }(v)=q_{\gamma }(u)$. Now, in the general case, (4.4) and (4.5) are not equivalent.

Evidently, the uniqueness with respect to (a) is stronger: the
lifting is unique along a given path iff it is unique along it through an arbitrary point above it.

{\bf Corollary 4.1.} The generated in accordance to (4.1) from an arbitrary transport along paths lifting of any path without self-intersection through any lying above it point is unique (in any one of the above-mentioned senses) along it.

{\bf Proof.} If $\gamma :J\to B$ is without intersections, then $\gamma (s)=\gamma (r), s,r\in J$ is equivalent to $r=s$, as a consequence of which for every $u\in \pi ^{-1}(\gamma (J))$ the set $q_{\gamma }(u)$ contains only one element, so the equalities $(4.3)-(4.5)$ are identically satisfied.\blacksquare

Due to corollary 4.1 there can be uniqueness of the mentioned lifts in a "middle" between local and global sense. For instance, we can say that the liftings generated by transports along paths are unique when they act on the class of paths without self-intersections.

Concluding the discussion of the uniqueness of lifting generated by transports along paths we shall prove

{\bf Proposition 4.3.} The lifting generated through (4.1) from I is unique
in a sense that for every path $\gamma :J\to B$, we have
\[
\bar{\gamma}_{\bar{\gamma}_{u,s^{(r),r}}}
=\bar{\gamma}_{u,s},\quad u\in \pi ^{-1}(\gamma (s)),\ r,s\in J,
\qquad (4.6)
\]
 i.e. the lift of any path through a point belonging to
its lifting, when the value of its parameter coincides with the one
describing this point, coincides with the initial lift of the path.

{\bf Proof.} (4.6) is equivalent to the basic property (2.2) of the
transports along paths: if $t\in J$, then
$\bar{\gamma}_{u,s}(t)=I^{\gamma }_{s\to t}(u)$ and
$\bar{\gamma}_{\bar{\gamma}_{u,s^{(r),r}}}(t)=I^{\gamma }_{r\to t}($
$_{u,s}(r))=I^{\gamma }_{r\to t}\circ I^{\gamma }_{s\to r}(u)$, so (4.6) and
(2.2) are equivalent.\blacksquare

{\bf Corollary 4.2.} If the lifting generated through (4.1) by I is globally
(resp. locally) unique (resp. along $\gamma )$ and $u,v\in \pi ^{-1}(\gamma
(J))$, then   $\bar{\gamma}_{u,r}(J)  $  $_{v,s}(J)=  $ or
$\bar{\gamma}_{u,r}(J)=\bar{\gamma}_{v,s}(J), r\in q_{\gamma }(u), s\in
q_{\gamma }(v)$ for every (resp. the given) path $\gamma $, i.e., in this
case two generated from I liftings of one and the same path do not intersect
with each other or coincide as sets.

{\bf Proof.} Let
$\bar{\gamma}_{u,r}(J)\cap\bar{\gamma}_{v,s}(J)\neq\emptyset$ and
$w\in\bar{\gamma}_{u,r}(J)\cap\bar{\gamma}_{v,s}(J)$. Then, due to the
global (resp.  local) uniqueness (resp. along $\gamma ) (4.5)$ is valid
because of proposition 4.2.  So
$\bar{\gamma}_{u,r}=\bar{\gamma}_{w,t}=\bar{\gamma}_{v,s}$ for $r\in
q_{\gamma }(u), s\in q_{\gamma }(v)$ and $t\in q_{\gamma }(w)$, i.e.
$\bar{\gamma}_{u,r}=\bar{\gamma}_{v,s}$, from where we get
$\bar{\gamma}_{u,r}(J)\bar{\gamma}=\bar{\gamma}_{v,s}(J).\blacksquare $

{\bf Corollary 4.3.} If the generated from I lifting is unique along a path
$\gamma , u\in \pi ^{-1}(\gamma (J)), v\in \pi ^{-1}(\pi (u))$ and $v\neq u$,
then   $\bar{\gamma}_{u,r}(J)\cap\bar{\gamma}_{v,s}(J)=\emptyset$, $r,s\in
q_{\gamma }(u)=q_{\gamma }(v)$, i.e., in this case two lifts of one and the
same path through two different points above it having equal projections do
not have a common point as sets.

{\bf Proof.} If we admit   $\bar{\gamma}_{u,r}(J)\cap
\bar{\gamma}_{v,s}(J)\neq\emptyset$ and
$w\in\bar{\gamma}_{u,r}(J)\cap\bar{\gamma}_{v,s}(J)$, then by proposition
4.2 (see (4.5)) $\bar{\gamma}_{u,r}=\bar{\gamma}_{w,t}=\bar{\gamma}_{v,s}$,
$r,s\in q_{\gamma }(u)=q_{\gamma }(v)$, $t\in q_{\gamma }(w)$. So, if we
choose $s^{*}\in q_{\gamma }(u)=q_{\gamma }(v)$ and use (4.3), we get $u=$
$\bar{\gamma}_{u,s^{*}}(s^{*})=\bar{\gamma}_{u,s}(s^{*})=
\bar{\gamma}_{v,s}(s^{*})=\bar{\gamma}_{v,s^{*}}(s^{*})=v$, which
contradicts $u\neq $v.\blacksquare

{\bf Proposition 4.4.} For every $s_{0}\in J$, for which the lifting
generated by I is unique along $\gamma $ through some $v_{0}\in \pi
^{-1}(\gamma (s_{0}))$, is fulfilled
\[
\pi ^{-1}(\gamma (J))
= \bigcap_{u\in \pi ^{-1}(\gamma(s_{0})) }\bar{\gamma}_{u,s}(J)
\qquad for\ \gamma (s)=\gamma (s_{0})\qquad (4.7)
\]
 for every $s\in J$ and if
$s=s_{0}$, then this equality is valid in spite of the above conditions.

{\bf Proof.} From the definition of the pointed uniqueness it follows that
$\bar{\gamma}_{v,s}=\bar{\gamma}_{v,s_{0}}$, $s,s_{0}\in q_{\gamma }(v_{0})$
for every $v\in \pi ^{-1}(\gamma (s_{0}))$ (see above). Then $v=$
$\bar{\gamma}_{v,s_{0}}(s_{0})=\bar{\gamma}_{v,s}(s_{0})$, from which we
find  $\bar{\gamma}_{v,s}=\bar{\gamma}_{\bar{\gamma}_{v,s_{0}}}$, so
\[
v\in
\bar{\gamma}_{v,s}(J) =\bar{\gamma}_{\bar{\gamma} _{v,s^{(s_{0}}}}(J)\subset
\bigcup_{u\in \pi ^{-1_{(}}} \bar{\gamma}_{(s_{0}}\gamma _{u,s}(J)
\]
and consequently
\[
 \pi ^{-1}(\gamma (J))
=\bigcup_{v\in \pi ^{-1} \gamma (J)) } v \subseteq
\bigcup_{u\in \pi ^{-1} \gamma(s_0) }\bar{\gamma}_{u,s}(J).
\]

Evidently, for $s=s_{0}$the last two formulae, as well as their proofs, are
identically valid and do not have any connection with the uniqueness of the
generated by I lifting.

 On the other hand, due to (4.2), we have the inclusion
\[
\bigcup_{u\in \pi ^{-1} \gamma(s_0)} \bar{\gamma}_{u,s}(J)
\subseteq \pi ^{-1}(\gamma (J)),
\]
from which, due to the previous inclusion, it follows (4.7).\blacksquare

{\bf Proposition 4.5.} If the generated by I lifting is unique along $\gamma
:J\to B$ and $E$ is a manifold, then for arbitrary fixed $s_{0},s\in J$ such
that $\gamma (s_{0})=\gamma (s)$, the family $\{\bar{\gamma}_{u,s}(J): u\in
\pi ^{-1}(\gamma (s_{0}))\}$ forms (one-dimensional) foliation of the
submanifold $\pi ^{-1}(\gamma (J))\subset $E.

{\bf Proof.} According to corollary 4.3, the sets from the considered family
do not have common points, so the proposition follows directly from the
definition of a foliation [11].\blacksquare

The above results show that the definition of a transport along paths in one
fibre bundle induces a structure of foliation in the set of liftings of an
arbitrary path in its base.

Now we shall consider the problem when in one (generally topological) fibre
bundle there exist transports along paths.

{\bf Theorem 4.1.} In the fibre bundle $(E,\pi ,B)$ there exist transports
along paths if and only if for every path $\gamma :J\to B$, every point $u\in
\pi ^{-1}(\gamma (J))$ and every $s\in q_{\gamma }(u)$ there exist a lift
$\bar{\gamma}_{u,s}:J\to E$ of $\gamma $ through $u$, for which
$\bar{\gamma}_{u,s}(s)=u$ and which is unique in a sense that for it (4.6)
is valid.

{\bf Proof.} The necessity follows directly from definition 4.1 and
propositions 4.1 and 4.3. On the opposite, let in $(E,\pi ,B)$ exist a lift
$l:(\gamma ,u,s)\to\bar{\gamma}_{u,s}$ of the paths in $B$ with the described
properties. For $s\in J$ and $u\in \pi ^{-1}(\gamma (s))$, we define the maps
$^{l}I^{\gamma }_{s\to t}:\pi ^{-1}(\gamma (s))\to \pi ^{-1}(\gamma (t))$ by
the equality
\[
^{l}I^{\gamma }_{s\to t}u
:=\bar{\gamma}_{u,s}(t),\quad s,t\in J,\ u\in \pi ^{-1}(\gamma (s)).
 \qquad (4.8)
\]

 On the one hand, if we put here $s=t$ and use
$\bar{\gamma}_{u,s}(s)=u$ and the arbitrariness of $u$, we get
$^{l}I^{\gamma }_{s\to s}=${\it id}$_{\pi ^{-1_{(\gamma (s))}}}$. On the
other hand, from the proof of proposition 4.3 it follows that (4.6) is
equivalent to $^{l}I^{\gamma }_{t\to r}\circ ^{l}I^{\gamma }_{s\to
t}=^{l}I^{\gamma }_{s\to r}$, so $^{l}I^{\gamma }_{s\to t}$is a transport
along $\gamma $ from $s$ to $t$, i.e., $^{l}I$ defined by $^{l}I:\gamma \to
^{l}I^{\gamma }:(s,t)\to ^{l}I^{\gamma }_{s\to t}$, is a transport along
paths in $(E,\pi ,B).\blacksquare $

 Theorem 4.1 has also a local variant expressed by

{\bf Theorem} ${\bf 4}{\bf .}{\bf 1}^\prime ${\bf .} In the fibre bundle
$(E,\pi ,B)$ along a given path $\gamma :J\to B$ there exists a transport
along it if and only if for every $u\in \pi ^{-1}(\gamma (J))$ and every
$s\in q_{\gamma }(u)$ there exists a lifting   $\bar{\gamma}_{u,s}:J\to E$
of $\gamma $ through $u$, for which   $\bar{\gamma}_{u,s}(s)=u$ and (4.6)
holds.

{\bf Proof.} The proof of  this theorem coincides with the proof of theorem
4.1 with the only difference that now the path $\gamma $ is
fixed.\blacksquare

The proof of theorem 4.1,  together with proposition 4.1, gives us a reason
to introduce

{\bf Definition 4.3.} If for  one lifting of paths equality (4.6) holds along
every (resp. a given) path $\gamma $, then for the defined by (4.8) transport
along paths we shall say that it is generated by this lifting (resp. along
$\gamma )$.

 {\bf Proposition 4.6.} If one transport along paths is generated through
(4.8) by some lifting of paths from $B$ to $E$, which is unique in a sense
that (4.6) is true, then the generated from this transport, in accordance
with definition 4.1, lifting of paths from $B$ to $E$ coincides with the
initial lifting  generating the considered transport along paths.

{ \bf Proof.} Let the transport along paths $^{l}I$ be generated by the lift
$l:\gamma \to \bar{\gamma}_{u,s}, u\in \pi ^{-1}(\gamma (J))$,
$_{u,s}(s)=u, s\in q_{\gamma }(u)$ in accordance with (4.8). Then, due to
definition 4.1 and proposition $4.1,{ } ^{l}I$ generates the lifting $^\prime
l:\gamma \to ^\prime $  $\bar{\gamma}_{u,s}$, so that
$^\prime\bar{\gamma}_{u,s}(t)=^{l}I^{\gamma }_{s\to t}$u. Comparing the last
equality with (4.8), we conclude that $^\prime\bar{\gamma}_{u,s}(t)=$
$\bar{\gamma}_{u,s}(t), t\in J$, i.e. $^\prime\bar{\gamma}_{u,s}=$
$\bar{\gamma}_{u,s}$ and consequently $^\prime l=$l.\blacksquare

{\bf Proposition 4.7.} If one lifting of paths from $B$ to $E$ is generated
from a transport along paths by (4.1), then this lifting generates in
accordance with $(4.8) a$ transport along paths coinciding with the transport
along paths generating the considered lifting.

{\bf Proof.} Let the lifting $\gamma \to\bar{\gamma}_{u,s}, u\in \pi
^{-1}(\gamma (J))$,   $\bar{\gamma}_{u,s}(s)=u, s\in q_{\gamma }(u)$ be
generated from the transport I, i.e.   $\bar{\gamma}_{u,s}(t)=I^{\gamma
}_{s\to t}$u.  Then, this lifting generates transport $^{l}I$ for which
$^{l}I^{\gamma }_{s\to t}u=\bar{\gamma}_{u,s}(t)=I^{\gamma }_{s\to t}u$, so
$^{l}I=$I.\blacksquare

\medskip
\medskip
 {\bf 5. CONNECTION GENERATED BY A TRANSPORT ALONG PATHS}

\medskip
In this section, the problem will be discussed of how one transport along paths, with the help of the generated by it lifting of paths, generates a connection in (part of) a given manifold.

Let $(E,\pi ,B)$ be a differential fibre bundle, $U$ be $k$-dimensional,
$1\le k\le \dim$, submanifold of $B$ and $U$ to be covered from a congruence
of paths $\gamma _{\lambda }:J\to U$ numbered with the $(k-1)$-dimensional
parameter $\lambda =(\lambda _{1},\ldots  ,\lambda _{k-1})\in \Lambda \subset
{\Bbb R}^{k}$, i.e.
$U=\cup_{\lambda\in\Lambda }{\bigl(}\gamma _{\lambda
}(J_{\lambda }){\bigr)}$
and
$\gamma _{\lambda }(J_{\lambda }) \cap \gamma_{\mu }(J_{\mu })= \emptyset$
for $\lambda ,\mu \in \Lambda , \lambda \neq \mu $.  Then, from (4.7) it
follows that for any fixed $x\in U$ is fulfilled
\[
\pi ^{-1}(U) =
\bigcup_{u\in \pi ^{-1} (x) } \bigl[
\bigcup_{\lambda\in\Lambda} \overline{(\gamma_\lambda)}_{u,s} (J_{\lambda
})\bigr],\quad s\in q_{\gamma }(u) \qquad (5.1)
\]
 where the sets in
 square brackets, which are $k$-dimensional manifolds, form $a k$-dimensional
foliation over $\pi ^{-1}(U) [11]$.

If $u\in \pi ^{-1}(U)$ and $T_{x}(M)$ is the tangent space to at $x\in M$ to
the manifold $M [3]$, then by definition [1,2,4] the vertical tangent space
to $\pi ^{-1}(U)$ at $u$ is
\[
T^{v}_{u}(\pi ^{-1}(U)):=T_{u}(\pi ^{-1}(\pi (u)))\subset T_{u}(\pi
^{-1}(U)).\qquad (5.2a)
\]

 Let:
\[
  T^{h-v}_{u}(\pi ^{-1}(U))
:=T_{u}{\bigl(} \bigcap_{\lambda \in\Lambda }
\overline{ (\gamma_\lambda) }_{u,s}
(J_{\lambda }){\bigr)}\subset T_{u}(\pi^{-1}(U)),\qquad (5.2b)
\]
\[
T^{0}_{u}(\pi ^{-1}(U))
:=T^{v}_{u}(\pi ^{-1}(U))\bigcap T^{h-v}_{u}(\pi ^{-1}(U))\qquad (5.2c)
\]
and
$T^{h}_{u}(\pi ^{-1}(U))$ be the direct complement of $T^{0}_{u}(\pi
^{-1}(U))$ in $T^{h-v}_{u}(\pi ^{-1}(U))$, i.e.
\[
T^{h-v}_{u}(\pi ^{-1}(U))=:T^{0}_{u}(\pi ^{-1}(U))\oplus T^{h}_{u}(\pi
^{-1}(U)).\qquad (5.2d)
\]

 {\bf Lemma 5.1.} For every $u\in \pi ^{-1}(U)$ is fulfilled
\[
T^{v}_{u}(\pi ^{-1}(U)) \bigcap T^{h}_{u}(\pi ^{-1}(U))
= \{{\bf 0}\}.\qquad (5.3)
\]
 where {\bf 0} is the zero element of $T_{u}(\pi ^{-1}(U))$.

{\bf Proof.} Let $u\in \pi ^{-1}(U)$ and
$\Omega :=T^{v}_{u}(\pi ^{-1}(U)) \cap T^{h}_{u}(\pi ^{-1}(U))$. If $\omega
\in \Omega $, then $\omega \in T^{h}_{u}(\pi ^{-1}(U))$ and according to
(5.2d), we have $\omega \in T^{h-v}_{u}(\pi ^{-1}(U))$. From here, taking
into account (5.2c) and $\omega \in T^{v}_{u}(\pi ^{-1}(U)) ($due to $\omega
\in \Omega )$, it follows $\omega \in T^{0}_{u}(\pi ^{-1}(U))$. This fact,
compared with $\omega \in T^{h}_{u}(\pi ^{-1}(U))$ shows that
$\omega \in T^{h}_{u}(\pi ^{-1}(U))\cap T^{0}_{u}(\pi ^{-1}(U))$. But, in
accordance with the definition for a direct sum, from (5.2d) it follows
$T^{h}_{u}(\pi ^{-1}(U))$  $T^{0}_{u}(\pi ^{-1}(U))=\{{\bf 0}\}$ and so
$\omega \in \{{\bf 0}\}$, i.e. $\omega ={\bf 0}$. Hence, if $\omega \in
\Omega $, then $\omega ={\bf 0}$, due to which $\Omega =\{{\bf 0}\}$, i.e.,
(5.3) is valid.\blacksquare

{\bf Proposition 5.2.} The subspaces $T^{h}_{u}(\pi ^{-1}(U))$ are
horizontal, i.e., they define a connection in the fibre bundle $(\pi
^{-1}(U),\pi \mid U,U)=(E,\pi ,B)\mid _{U}$, iff
\[
 T_{u}(\pi ^{-1}(U))
=T^{v}_{u}(\pi ^{-1}(U))\oplus T^{h}_{u}(\pi ^{-1}(U)).\qquad (5.4)
\]

 {\bf Remark.} In the case of locally trivial fibre
bundles an elementary check (see (5.2) and the dimensions of the defined by
it objects) shows the equivalence of (5.4) and the condition
\[
 T^{0}_{u}(\pi ^{-1}(U))=\{0\},\qquad (5.4^\prime )
\]

 {\bf Proof.} The proposition is a direct consequence of the definitions
(5.2), lemma 5.1 and the definition of a connection in arbitrary fibre
bundles $[2-4].\blacksquare $

A priori it is clear that the validity of the condition (5.4) depends on the used transport along paths I, which are used to construct liftings of paths appearing above, as well as on the choice of the family of paths $\{\gamma _{\lambda }\}$.

In particular, if  $U=B$,  and  (5.4)  is  fulfilled,  then  the $\dim(B)$-dimensional distribution $u\to T^{h}_{u}(E)=T^{h}_{u}(\pi ^{-1}(B)), u\in E$ defines a connection in $(E,\pi ,B)$ and it is almost evident that the defined  by it parallel transport (see [3,4]) along the paths $\gamma _{\lambda }, \lambda \in \Lambda $ coincides with the initial transport I along them.

\medskip
 {\bf 6.CONCLUSION}

\medskip
In this work, we have defined and investigated some properties of the transports along paths in general fibre bundles. As examples of such transports we pointed out the linear transports along paths in vector bundles [6] and, in particular as their special case, the generated by derivations of tensor algebras over a manifold transports along paths [5]. Here, we were not concern with the ties of the theory developed with the ones of connections, parallel transports, bundle morphisms etc., items which will be studied elsewhere.

\medskip
\medskip
 {\bf ACKNOWLEDGEMENTS}

\medskip
The author expresses his gratitude to Prof. Vl. Aleksandrov (Institute of Mathematics of Bulgarian  Academy of Sciences) for constant interest in this work and stimulating discussions.

This research was partially supported by the Fund for Scientific Research of Bulgaria under contract Grant No. $F 103$.

\medskip
\medskip
 {\bf REFERENCES}

\medskip
1.  Greub W., S. Halperin, R. Vanstone, Connections, Curvature, and Cohomology, vol.1, vol.2, Academic Press, New York and London, $1972, 1973$.\par
2.  Dubrovin B.A., S.P. Novikov, A.T. Fomenko, Modern geometry, Nauka, Moscow, 1979 (In Russian).\par
3.  Kobayashi S., K. Nomizu, Foundations of differential geometry, vol.1, Interscience publishers, New-York-London, 1963.\par
4.  Warner F.W., Foundations of differentiable manifolds and Lie groups, Springer Verlag, New York-Berlin-Heidelberg-Tokyo, 1983.\par
5.  Iliev B.Z., Parallel transports in tensor spaces generated by derivations of tensor  algebras, Communication JINR, $E5-93-1$, Dubna, 1993.\par
6.  Iliev B.Z., Linear transports along paths in vector bundles. I. General theory, Communication JINR, $E5-93-239$, Dubna, 1993.\par
7.  Steenrod N., The topology of fibre bundles, 9-th ed., Princeton Univ. Press, Princeton, $1974 (1-st ed. 1951)$.\par
8.  Viro O.Ya., D.B. Fuks, I. Introduction to homotopy theory, In: Reviews of science and technic, sec. Modern problems in mathematics. Fundamental directions, vol.24, Topology-2, VINITI,
Moscow, $1988, 6-121 ($In Russian).\par
9.  Husemoller D., Fibre bundles, McGrow-Hill Book Co., New York-St. Louis-San Francisco-Toronto-London-Sydney, 1966.\par
10.  Iliev B.Z., Linear transports along paths in vector bundles. II. Some applications, Communication JINR, $E5-93-260$, Dubna, 1993.\par
11.  Tamura I, Topology of foliations, Mir, Moscow, 1979 (In Russian; translation from Japanese).\par
$12. Hu$ Sze-Tsen, Homotopy Theory, Academic Press, New York-London, 1959.\par
13.  Iliev B. Z., Consistency between metrics and linear transports along curves, Communication JINR, $E5-92-486$, Dubna, 1992.

\newpage

\noindent Iliev B. Z.\\[5ex]

\noindent Transports along Paths in Fibre Bundles\\
 I. General Theory\\[5ex]

\medskip
\medskip
Transports along path in fibre bundles are axiomatically introduced. Their
general functional form and some their simple properties are investigated.
The relationships of the transports along paths and lifting of paths are
studied.\\[5ex]

\medskip
 The investigation has been performed at the Laboratory of Theoretical
Physics, JINR.

\end{document}